\documentclass[12pt,a4paper]{article}

\usepackage{amsfonts, amsmath, amsthm}
\usepackage{graphicx, float}

\newtheorem{theorem}{Theorem}
\newtheorem{proposition}[theorem]{Proposition}
\newtheorem{lemma}[theorem]{Lemma}
\newtheorem{corollary}[theorem]{Corollary}

\theoremstyle{remark}
\newtheorem{example}{Example}

\newcommand{\R}{\mathbb{R}}
\newcommand{\bl}{\boldsymbol{\ell}}
\newcommand{\ind}{\mathbb{I}}
\newcommand{\var}{\operatorname{Var}}

\title{Generalized R\'enyi statistics}

\author{P\'eter Kevei\footnote{Bolyai Institute, University of Szeged,
Hungary.
E-mail: \texttt{kevei@math.u-szeged.hu}} \and
L\'aszl\'o Viharos\footnote{Bolyai Institute, University of Szeged, Hungary.
E-mail: \texttt{viharos@math.u-szeged.hu}}}

\date{}

\begin{document}

\maketitle

\begin{abstract}
In R\'enyi's representation for exponential order statistics,
we replace the iid exponential sequence with any iid sequence, and 
call the resulting order statistic \emph{generalized R\'enyi statistic}. 
We prove that by randomly reordering the variables in the generalized
R\'enyi statistic, we obtain in the limit a sequence of iid exponentials.
This result allows us to propose a new model for heavy-tailed data.
Although the new model is very close to the classical iid framework,
we establish that the Hill estimator is weakly consistent and 
asymptotically normal without any further assumptions on the 
underlying distribution or on the number of upper order statistics
used in the estimator.


\end{abstract}

\section{Introduction and main results}

Let $Z_1, Z_2, \ldots$ be independent identically distributed (iid)
random variables with finite positive mean $\gamma = E (Z_1)$.
Here we are interested in the 
\emph{generalized R\'enyi statistics}, defined as
\begin{equation} \label{eq:generalized}
X_{k,n} =\sum_{j=1}^k\frac{Z_j}{n+1-j}, \quad k=1,2,\ldots,n.
\end{equation}
If the $Z$'s have a common exponential distribution then it is well-known 
from R\'enyi's representation \cite{Renyi}
that the vector $(X_{1,n}, X_{2,n}, \ldots, X_{n,n})$ has the same distribution 
as the order statistics of $n$ iid exponentials.
R\'enyi's representation theorem has numerous applications in the 
theory of records and order statistics, see 
e.g.~Nevzorov \cite[Representation 3.4]{Nevzorov}, 
Shorack and Wellner \cite[Chapter 21]{Shorack}, or 
de Haan and Ferreira \cite[Section 2.1]{deHaan}.

Consider a uniform random permutation $(\delta_1,\ldots,\delta_n)$
of $(1,\ldots, n)$, 
independent of $X_{1,n},  \ldots,  X_{n,n}$ in \eqref{eq:generalized}. Then
$X_{\delta_1,n},\ldots, X_{\delta_n,n}$ are identically distributed.
Moreover, from R\'enyi's representation we see that 
$X_{\delta_1,n},\ldots, X_{\delta_n,n}$ are iid exponentials
if and only if $Z_1, Z_2, \ldots$ are exponentials.
In our main results, we establish that this remains true asymptotically in
the more general model \eqref{eq:generalized}, demonstrating
the universality property of the exponential distribution.
Throughout, $\stackrel{\mathcal{D}}{\longrightarrow}$ denotes
convergence in distribution.

\begin{theorem} \label{thm:exponential}
Assume that  $Z_1, Z_2, \ldots$ in \eqref{eq:generalized} are 
iid random variables such  that $E (Z_1 ) = \gamma > 0$, and 
let $(\delta_1, \ldots, \delta_n)$ be a uniform random permutation
of $(1,\ldots, n)$ independent of the $Z$'s. Then for any $k \geq 1$
\[
\left( X_{\delta_1, n}, \ldots, X_{\delta_k,n} \right) 
\stackrel{\mathcal{D}}{\longrightarrow}
(Y_1, \ldots, Y_k), \quad \text{as } \, n \to \infty,
\]
where $Y_1, \ldots, Y_k$ are independent exponential random variables 
with mean $\gamma$.
\end{theorem}

Under stronger moment conditions, the moments also converge.

\begin{theorem} \label{thm:moments}
Assume that $Z_1,Z_2,\ldots$ are iid random variables with mean $\gamma$ and 
$E(|Z_1|^t)<\infty$ for all $t>0$.
Then $E(X_{\delta_1,n}^k)\to k!\gamma^k$, for any $k=1,2,\dots$, and 
$X_{\delta_1,n},\ldots,X_{\delta_n,n}$ are pairwise asymptotically uncorrelated,
that is $E (X_{\delta_1, n} X_{\delta_2,n} ) \to \gamma^2$, as $n \to \infty$.
\end{theorem}

\smallskip

Based on these distributional limit results, we propose a new model for 
heavy-tailed data. We assume that the order statistics 
$W_{1,n} \le\cdots \le W_{n,n}$ of heavy-tailed data can be represented in 
distribution as
\begin{equation} \label{eq:modelheavy}
\begin{split}
& (W_{k,n}:k=1,\ldots,n) \stackrel{\mathcal{D}}{=} (C e^{X_{k,n}}: k=1,\ldots,n), \ 
\text{ where $C$ is a positive}  \\
& \text{scale factor and $X_{k,n}$'s are given in 
\eqref{eq:generalized} with $Z_1,Z_2,\ldots$ iid} \\
& \text{\emph{nonnegative} random 
variables having finite mean $\gamma = E (Z_1) > 0$}.
\end{split}
\end{equation}

Note that we do not assume a priori that the $Z$'s in \eqref{eq:generalized}
are nonnegative. In particular, Theorems \ref{thm:exponential} and 
\ref{thm:moments} hold for general $Z$. However, in model \eqref{eq:modelheavy} 
the $Z$'s must be nonnegative to ensure that the $X_{k,n}$'s are ordered.

An equivalent formulation of the assumptions in \eqref{eq:modelheavy} is that 
the scaled log-spacings, $(\log W_{k,n} - \log W_{k-1,n}) (n-k+1)$, 
for $k = 1,2,\ldots, n$, are iid with  a finite expectation $\gamma$,
where $W_{0,n} = C$.

\smallskip 

We demonstrate that the new model includes the iid strict Pareto model
which was first proposed by Pareto \cite{Pareto} to describe various income distributions.
Let $Z_1, Z_2, \ldots$ be iid exponentials with mean $\gamma > 0$. 
By R\'enyi's representation, the sequence 
$X_{1,n} < \ldots < X_{n,n}$ has the same 
distribution as ordered iid exponentials with mean $\gamma$.
Furthermore,  
\[
P( C e^{Z_1} > x) = \left( \frac{C}{x} \right)^{1/\gamma}, \quad x \geq C,
\]
i.e.~$C e^{Z_1}$ has a strict Pareto distribution. 
We note that there is a similar connection between the \emph{general 
exponential model} and the class of \emph{regularly varying distributions},
see Cs\"{o}rg\H{o} and Mason \cite[Fact 1.1]{CsM}.
Therefore, $(W_{k,n}: k= 1,\ldots, n)$ has the same distribution as the order statistics 
of an iid strict Pareto sample. This particular example also shows that 
the parameter $1/\gamma =1/E(Z_1)$ plays the role of the tail index in the proposed model.
Furthermore, Theorem \ref{thm:exponential} implies that this is the only example 
satisfying the assumptions of the iid model and model \eqref{eq:modelheavy}.

In the iid model with regularly varying tail the scaled log-spacings can 
be approximated by iid exponentials under suitable assumptions. Results 
in this direction were obtained by Beirlant et al.~\cite{Beirlant}.

\smallskip 

The estimation of the tail index in heavy-tailed models 
is crucial for many applications. 
The popular Hill estimator for $\gamma$ was introduced by Hill 
\cite{Hill} in the iid setup about 50 years ago. 
Despite its long history, the study of the 
asymptotic properties of the Hill estimator and its variants remains
an active area of research, see e.g.~Paulauskas and Vai\v{c}iulis \cite{Paulauskas}, 
Einmahl and He \cite{Einmahl}, or Wang and Resnick \cite{WangResnick}.

The asymptotic normality of the Hill estimator, 
centered by the true parameter,
holds under a second-order regular variation condition, see 
Resnick \cite[Section 9.1.2]{Resnick}, and the references therein.
Second-order regular variation is a rather strict assumption on the tail behavior 
of the distribution, and cannot be verified in practice. Furthermore, a key challenge 
in the application  of the Hill estimator is the choice of the 
parameter $k_n$, which determines
the number of upper order statistics included in the estimation. In 
model \eqref{eq:modelheavy} both issues are resolved. The Hill estimator 
simplifies to the average of iid random variables.
As a result, consistency and asymptotic normality hold without any 
further assumptions. Thus, while the proposed model is, by Theorem 
\ref{thm:exponential}, very close to the iid framework, the properties of 
the corresponding Hill estimator are much better. 
\smallskip

In Section \ref{sect:stat} we examine statistical inference for the 
new model. All the proofs are presented in Section \ref{sect:proofs}.

\section{Statistical inference} \label{sect:stat}

We investigate different estimators of the crucial parameter $\gamma$.

\subsection{Quantile estimator}

Based on the order statistics $X_{1,n}, \ldots, X_{n,n}$, define the empirical quantile 
functions $Q_n^{(1)}(s)$ and $Q_n^{(2)}(s)$ as
\begin{equation} \label{eq:quant}
Q_n^{(1)}(s) :=X_{\lceil ns \rceil,n},\quad Q_n^{(2)}(s) 
:=e^{Q_n^{(1)}(s)}, \qquad s \in (0,1),
\end{equation}
where $\lceil \cdot \rceil$ stands for the upper integer part.

The empirical quantiles allow us to introduce a natural estimate for $\gamma$. 
Define for $s \in (0,1)$
\begin{equation}
\widetilde \gamma_n(s) = \frac{Q_n^{(1)}(s)}{-\log (1-s)}.
\end{equation}
We prove that $\widetilde \gamma_n(s)$ is a weakly consistent and asymptotically 
normal estimator of $\gamma$.

\begin{theorem} \label{thm:Qn}
Assume that \eqref{eq:generalized} holds with nonnegative 
iid random variables $Z_1, Z_2, \ldots$ such that 
$\var (Z_1) = \sigma^2 < \infty$. Then as $n \to \infty$
\[
\sqrt{n} (\widetilde \gamma_n(s) - \gamma) 
\stackrel{\mathcal{D}}{\longrightarrow}
N(0, \sigma^2 h(s))
\]
for any $s \in (0,1)$, where 
\begin{equation} \label{eq:sigma}
h(s) =  \frac{{s}}{(1-s) (\log (1-s))^2}.
\end{equation}
\end{theorem}

As a consequence, a $1-\varepsilon$ asymptotic confidence interval is
\[
\left[
\widetilde \gamma_n(s) -\frac{\widetilde \sigma_n \sqrt{h(s)} x_\varepsilon}{\sqrt{n}},\ 
\widetilde \gamma_n(s) +\frac{\widetilde \sigma_n \sqrt{h(s)} x_\varepsilon}{\sqrt{n}}
\right],
\]
where $\widetilde \sigma_n$ is a consistent estimator of $\sigma$, 
and $x_\varepsilon$ is the $1 -\frac{\varepsilon}{2}$-th quantile of the 
standard normal distribution.

Asymptotic normality implies weak consistency, thus 
\[
Q_n^{(1)}(s)\stackrel{\mathbb{P}}{\longrightarrow} -\gamma\log(1-s),
\]
where $\stackrel{\mathbb{P}}{\longrightarrow}$
stands for convergence in probability. It means that the limiting quantile 
function of the generalized R\'enyi 
statistics is the quantile function of the exponential distribution.
Here and later on, any nonspecified limit relation is meant as $n \to \infty$.
Thus, by \eqref{eq:quant}
\[
\frac{Q_n^{(2)}(1-st)}{Q_n^{(2)}(1-s)}
\stackrel{\mathbb{P}}{\longrightarrow}  t^{-\gamma},
\]
for any $0<s,st<1$. The latter convergence reflects the regular 
variation property.

\begin{figure}
\begin{center}
\includegraphics[width=0.49\textwidth]{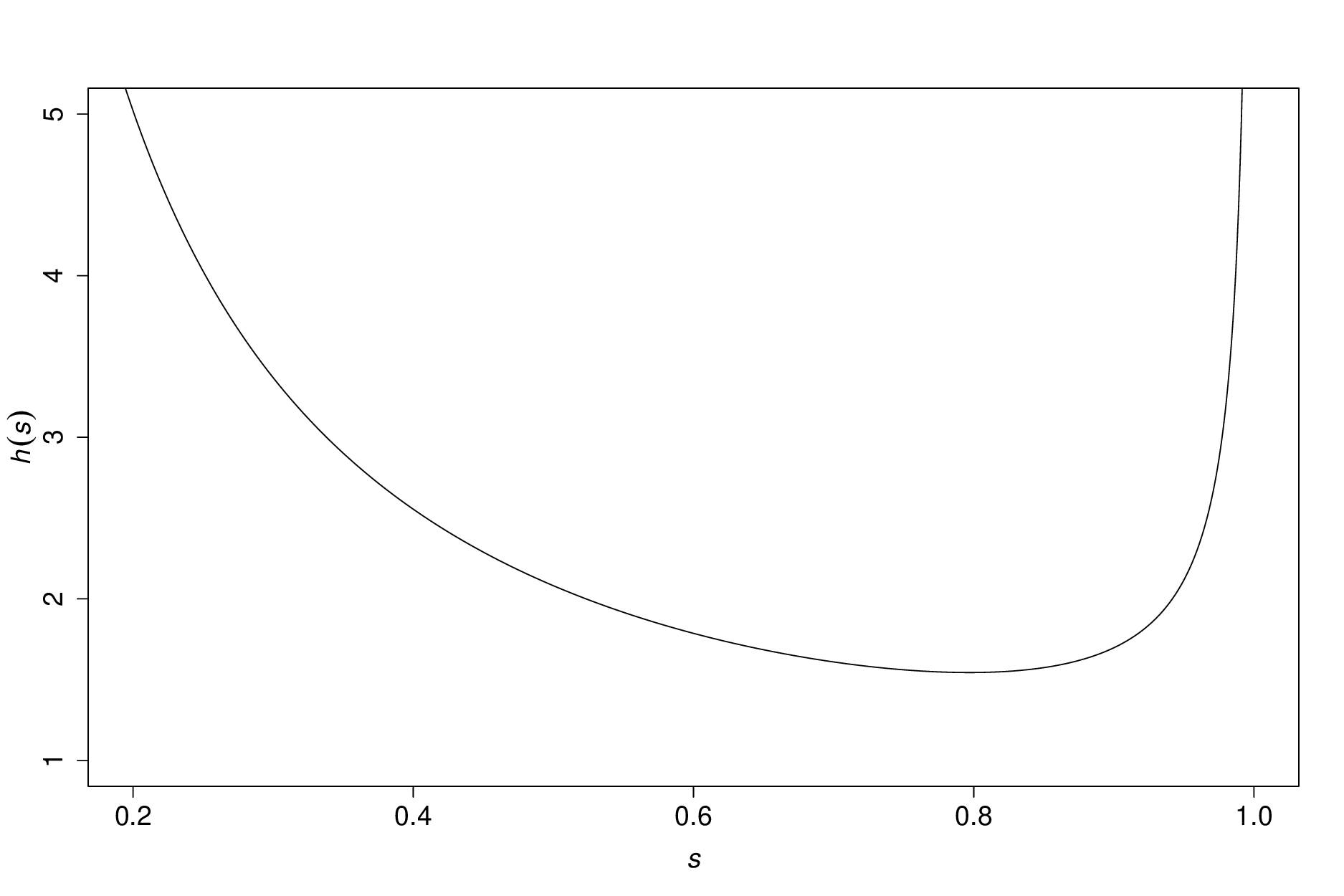}
\hspace*{-0.2cm}
\includegraphics[width=0.49\textwidth]{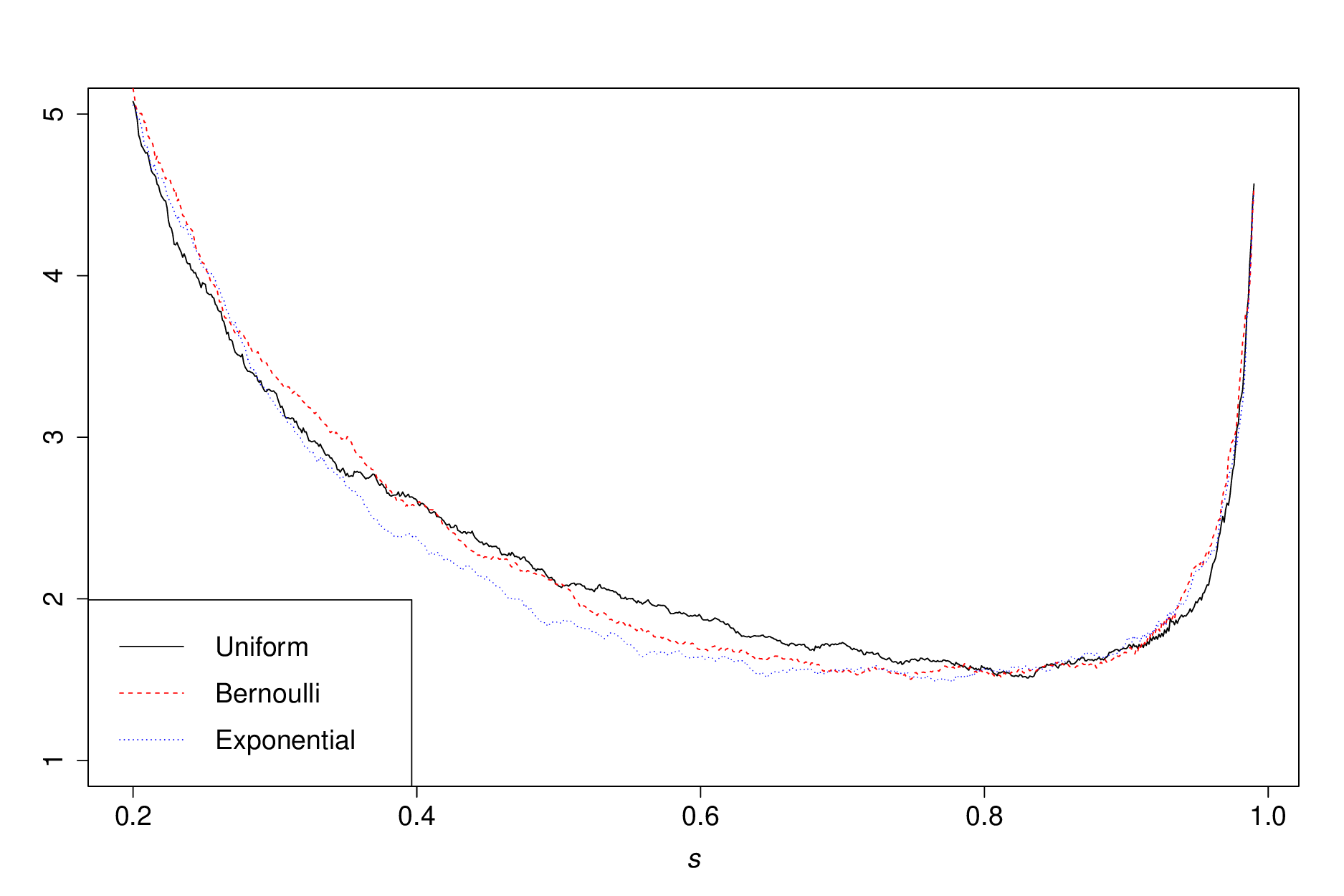}
\end{center}
\vspace*{-0.4cm}
\caption{
Left: The function $h(s)$ appearing in the asymptotic variance of 
$\widetilde \gamma_n(s)$.
Right: Empirical variances of $\sqrt{n} (\widetilde \gamma_n(s) - \gamma) / \sqrt{ \var(Z)}$, 
with $n= 1000$, for 1000 independent simulations. The common distribution of the $Z$'s 
in \eqref{eq:generalized} is uniform, Bernoulli, and exponential, respectively.}
\label{fig:sigma}
\end{figure}

In Figure \ref{fig:sigma} we plot the function $h(s)$ 
that appears in the asymptotic  variance. A short calculation gives 
that $h(s)$ has a unique minimum, $h(s_0) \approx 1.544$, 
attained at $s_0 \approx 0.797$, 
which is the 
unique solution on $(0,1)$ to the equation $\log \frac{1}{1-s} = 2 s$.

We simulated the estimator $\widetilde \gamma(s)$ for $s \in (0.2, 0.99)$ 
with step size $0.001$ in $s$, and $Z$'s having Uniform$(0,2\gamma)$, 
Exponential$(1/\gamma)$, and Bernoulli$(\gamma)$ distributions, 
each with $\gamma = 0.5$. 
We simulated $1000$ times with sample size $n = 1000$.
In Figure \ref{fig:sigma}, we display the simulated asymptotic variances 
normalized by $\var (Z)$, that is the empirical variance of 
$\sqrt{n} (\widetilde \gamma_n(s) - \gamma) / \sqrt{ \var(Z)}$, 
for 1000 independent simulations. 
Even with relatively small sample size, the 
simulated variances are very close to the theoretical values.

\subsection{Maximum likelihood estimator} \label{ssect:maxl}

Next, we assume that the common distribution of $Z$'s has density function $g$.
In this case the likelihood function can be determined explicitly as follows.

\begin{proposition} \label{prop:likelihood}
Assume that \eqref{eq:modelheavy} holds, and the random variables $Z_1,Z_2,\ldots$ 
have density $g$.
Then the conditional distribution of $(W_{n-k+1,n},\ldots,W_{n,n})$ given 
$W_{n-k,n} = w_{n-k}$ is absolutely continuous with density function
\begin{equation} \label{eq:likelihood}
\begin{split}
& h(w_{n-k+1}, w_{n-k+2},\ldots, w_n| w_{n-k}) \\
& =k!\prod_{j=n-k+1}^ng\big((n-j+1)(\log w_j -\log w_{j-1})\big) w_j^{-1},
\end{split}
\end{equation}
for $C\le w_{n-k}\le w_{n-k+1} \le \cdots \le w_n$, $k=1,\ldots,n$, $w_0 =C$.
\end{proposition}

In applications, we can use this likelihood function to fit parametric models to 
the scaled log-spacings of the original data. 
In particular, we can construct estimates for $\gamma$ in 
model \eqref{eq:modelheavy}.
We study the following parametric distributions.

\begin{example}
If the random variables $Z_k$'s are exponential with density
$g(x) =\frac{1}{\gamma}e^{-\frac{x}{\gamma}}$, $x\ge0$, then
the likelihood \eqref{eq:likelihood} is maximized for
\[
\widehat\gamma^{ML}_n =\widehat\gamma_n(k)=\frac{1}{k}\sum_{j=1}^k \log W_{n+1-j,n} -\log 
W_{n-k, n},
\]
which is exactly the Hill estimator in \eqref{eq:def-hill-W}.
\end{example}

\begin{example} \label{ex:gamma}
Assume that 
\begin{equation} \label{eq:gamma}
Z_k \text{ has gamma distribution with parameters $r$ and $r/\gamma$}, 
\ r > 0, \gamma > 0.
\end{equation}
Then the maximum likelihood estimator again coincides with the Hill estimator, 
$\widehat\gamma^{ML}_n=\widehat\gamma_n(k)$.
\end{example}

\begin{example}
If $Z_k$'s are uniform on $(0, 2\gamma)$, that is 
$g(x) =\frac{1}{2\gamma}$, $0<x<2\gamma$, then
the maximum likelihood estimator of $\gamma$ is
\[
\widehat\gamma^{ML}_n =\frac{1}{2}\max_{\,n-k_n+1\le j\le n}(n-j+1)(\log W_{j,n} -\log 
W_{j-1,n}).
\]
\end{example}

\subsection{Hill's estimator} \label{sect:asnorm}

The Hill estimator is the most widely studied estimator of the tail index. 
In the iid case, 
asymptotic normality of the Hill estimator, along with its various generalizations,
has been established under different assumptions by Hall \cite{Hall}, 
Segers \cite{Segers}, Gomes and Martins \cite{GomesMartins}, and 
Paulauskas and Vai\v{c}iulis \cite{Paulauskas}, to name a few. 
For a textbook treatment, see 
\cite[Theorem 3.2.5]{deHaan},
\cite[Section 9.1.2]{Resnick}, and 
the references therein. Essentially, asymptotic normality with centering 
with the true parameter $\gamma$, holds under a second-order 
regular variation condition.
In the classical iid case, the Hill estimator is calculated from $k$
upper order statistics, which are dependent random variables, complicating 
the analysis.

In our model \eqref{eq:modelheavy}, the Hill estimator takes the form
\begin{equation} \label{eq:def-hill-W}
\widehat\gamma^H_n(k)
=\frac{1}{k}\sum_{j=1}^{k}j(\log W_{n-j+1,n} -\log W_{n-j,n})
\stackrel{\mathcal{D}}{=}\frac{1}{k}\sum_{j=1}^{k}Z_{n-j+1},
\end{equation}
which is simply the average of iid random variables. As a result, 
Theorems \ref{thm:normal}
and \ref{thm:large} follow immediately from standard theory. 
In particular, asymptotic normality holds without any further assumptions 
on the underlying distribution or on the subsequence $k = k_n$.

\begin{theorem} \label{thm:normal}
Assume that \eqref{eq:modelheavy} holds.
If $\var (Z_1) = \sigma^2 < \infty$, $1 \leq k_n \leq n$, 
and $k_n \to \infty$, then
\[
\sqrt{k_n}\, \big(\widehat\gamma^H_n(k_n) -\gamma\big) 
\stackrel{\mathcal{D}}{\longrightarrow} \,{N}(0,\sigma^2) \quad \text{as } 
n \to \infty.
\]
\end{theorem}

As a consequence, in model \eqref{eq:modelheavy} 
a $1-\varepsilon$ asymptotic confidence interval is
\begin{equation} \label{eq:confint-new}
\left[
\widehat\gamma^H_n(k_n) -\frac{\widehat \sigma_{k_n} x_\varepsilon}{\sqrt{k_n}},\ 
\widehat\gamma^H_n(k_n) +\frac{\widehat \sigma_{k_n} x_\varepsilon}{\sqrt{k_n}}
\right],
\end{equation}
where the estimate $\widehat \sigma_{k_n}^2$ stands for the empirical variance of 
the scaled log-spacings $(n-j+1)(\log W_{j,n} - \log W_{j-1,n})$,
$j=1,2,\ldots, k_n$, and $x_\varepsilon$ is the $(1 -\frac{\varepsilon}{2})$-th 
quantile of the standard normal distribution.

In the iid setup, under appropriate assumptions 
(see e.g.~\cite[Theorem 2]{Hall}, \cite[Theorem 3.2.5]{deHaan}, \cite[Proposition 9.3]{Resnick})
\begin{equation} \label{eq:confint-Hill}
\left[
\widehat\gamma^H_n(k_n) -\frac{\widehat\gamma^H_n(k_n)x_\varepsilon}{\sqrt{k_n}},\ 
\widehat\gamma^H_n(k_n) +\frac{\widehat\gamma^H_n(k_n)x_\varepsilon}{\sqrt{k_n}}
\right],
\end{equation}
is a $1-\varepsilon$ asymptotic confidence interval for $\gamma$.
We observe that the length of the confidence interval cannot be controlled by model 
parameters other than $\gamma$. However, if $\sigma \leq \gamma$ then in 
model \eqref{eq:modelheavy}, the resulting asymptotic confidence interval 
is narrower compared to the classical iid setup.

\begin{figure}
\begin{center}
\includegraphics[width=0.49\textwidth]{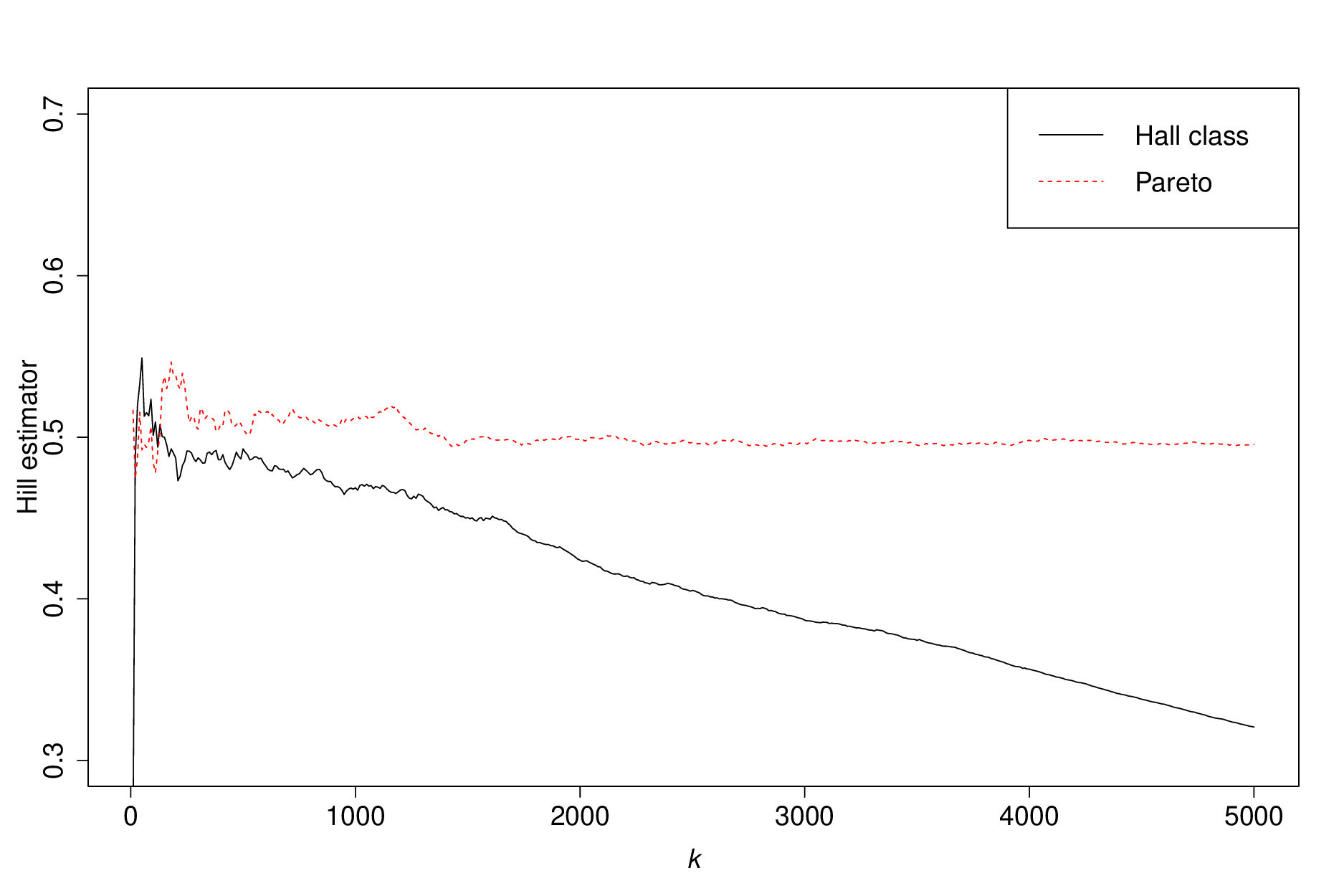}
\hspace*{-0.2cm}
\includegraphics[width=0.49\textwidth]{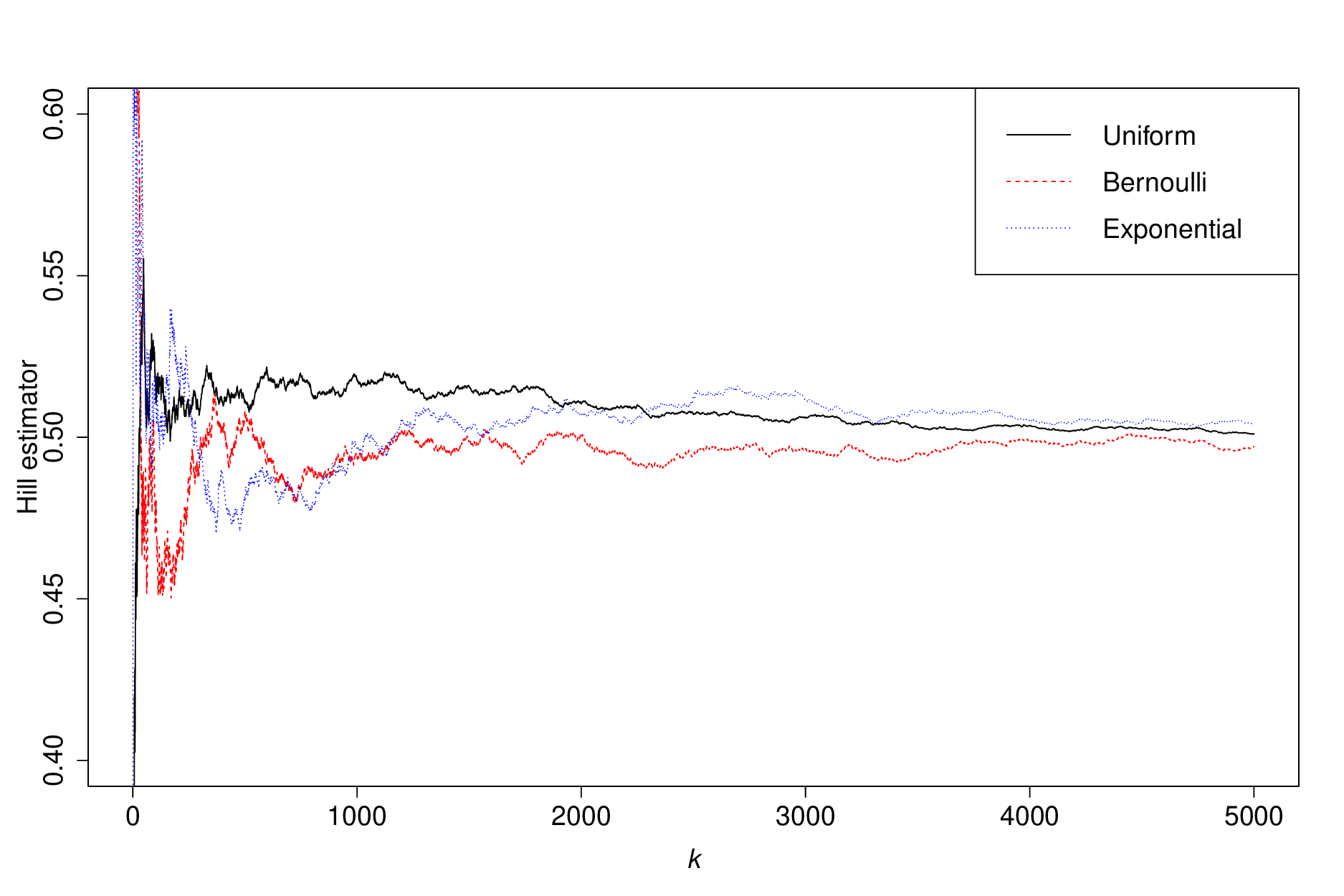}
\end{center}
\vspace*{-0.4cm}
\caption{The Hill estimator in the classical iid model (left) and in 
model \eqref{eq:modelheavy} (right).}
\label{fig:hill}
\end{figure}

In a simulation study, we compared the performance of the Hill estimator in the 
classical iid model and in model \eqref{eq:modelheavy}. On the left side  of 
Figure \ref{fig:hill}, we show the graph of the Hill estimator 
$\{(k, \widehat \gamma^H_n(k)) : k= 1, \ldots, n \}$, 
for $n=5000$ iid random variables, using two distributions: 
(i) strict Pareto distribution with $\gamma = \tfrac{1}{2}$, 
and (ii) a small perturbation of the Pareto, 
a distribution with quantile function
$Q(1-u) = u^{-{1}/{2}} ( 1 + \tfrac{1}{2} u)$. The latter distribution
belongs to the Hall class in \cite{Hall}.
We note again that if we choose a Pareto distribution in the iid model, 
the resulting model also satisfies the conditions of model \eqref{eq:modelheavy}. 
On the right side of Figure \ref{fig:hill}, we show the 
graph of the Hill estimator for model \eqref{eq:modelheavy},
using three distributions for $Z$: uniform, Bernoulli, and exponential, each 
with mean $\gamma = \tfrac{1}{2}$. In all cases $n=5000$. 
It is apparent from the figures that the Hill estimator performs better in model \eqref{eq:modelheavy}.
We observe that a larger $k_n$ leads to better estimates, with $k_n =n$ 
being the optimal choice,
which is theoretically clear from Theorem \ref{thm:normal}.

\begin{figure}
\begin{center}
\includegraphics[width=0.49\textwidth]{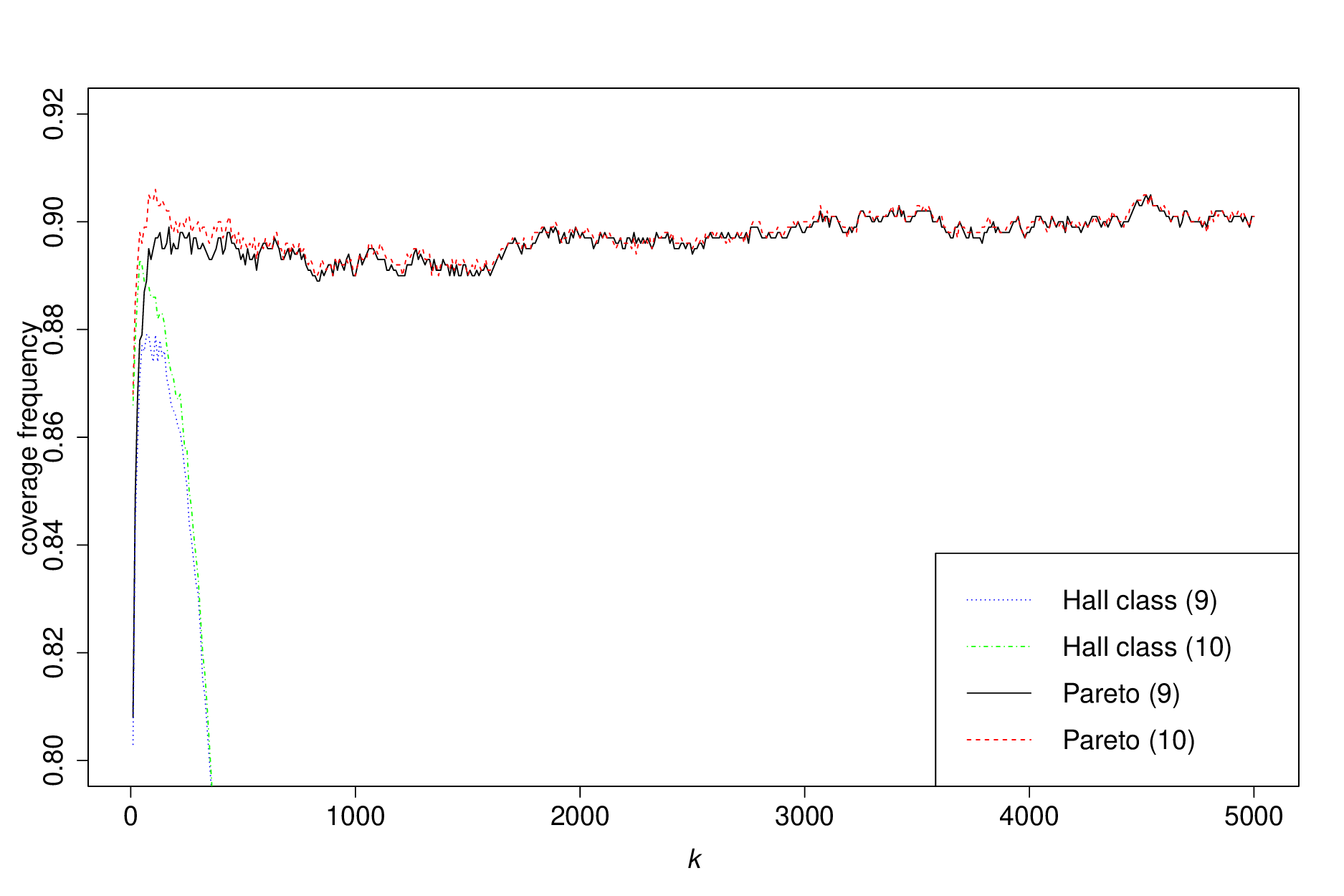}
\hspace*{-0.2cm}
\includegraphics[width=0.49\textwidth]{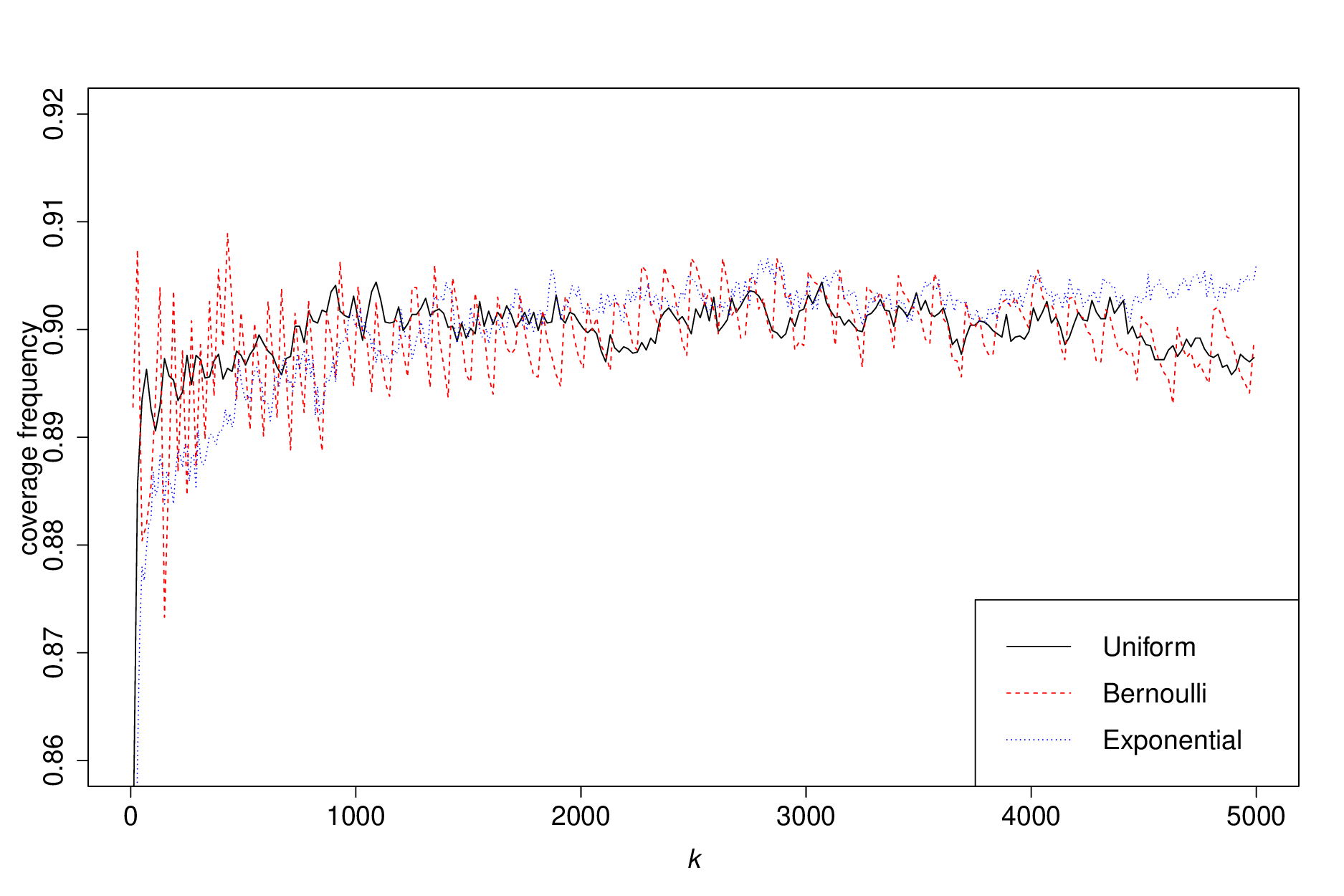}
\end{center}
\vspace*{-0.4cm}
\caption{Coverage frequencies of the asymptotic 
confidence intervals for $\varepsilon = 0.1$ in the classical iid model (left) and in model \eqref{eq:modelheavy} (right). In the iid model we used confidence interval 
\eqref{eq:confint-new} and \eqref{eq:confint-Hill}.}
\label{fig:confint}
\end{figure}

We repeated the simulation $10,000$ times and calculated the coverage 
frequencies of the asymptotic confidence intervals in \eqref{eq:confint-new}
and \eqref{eq:confint-Hill} with $\varepsilon = 0.1$. 
In the classical iid case, 
we calculated both with interval \eqref{eq:confint-new} and \eqref{eq:confint-Hill}.
The resulting coverage frequencies are plotted in Figure \ref{fig:confint}. It 
is clear from Figures \ref{fig:hill} and \ref{fig:confint}  that, 
in the Hall class, the Hill estimator works only for small $k$.

We note that in the classical iid setup we obtained similar coverage frequencies 
with both asymptotic confidence intervals  \eqref{eq:confint-new} and 
\eqref{eq:confint-Hill}. This is not surprising, since in the strict Pareto 
case $\widehat \sigma_{k_n}$ is a consistent estimator of $\gamma$.
We expect that $\widehat \sigma_{k_n}$ is a consistent 
estimator of $\gamma$ under suitable conditions on the slowly varying function 
in the tail. A proof could be based on the uniform approximation of 
the scaled log-spacings in a regularly varying iid model obtained in 
\cite{Beirlant}. This problem remains for further research.

Next, we turn to large deviations.
Consider the moment generating function
\[
M(t) = E\big(e^{t Z_1}\big),\quad-\infty< t<\infty,
\]
and define the rate function
\[
I(z) =\sup_{t\in\mathbb R}\big(zt -\log M(t)\big),\quad-\infty<z<\infty.
\]
The following theorem is a reformulation of Cram\'er's large deviation theorem for 
$\widehat\gamma_n(k_n)$, see e.g.~den Hollander \cite[Section I.3-4]{denHollander}.

\begin{theorem} \label{thm:large}
Assume the conditions of Theorem \ref{thm:normal}.
Then, if $M$ is finite in an open neighborhood of zero, we have
\[
\lim_{n\to\infty} \frac{1}{k_n}\log P\left(\widehat\gamma_n(k_n)\ge y\right)
= -\inf_{x\ge y} I(x),\quad-\infty<y<\infty.
\]
\end{theorem}

In special cases, we can calculate  the rate function explicitly.

\begin{corollary} \label{cor:ldp}
Assume that \eqref{eq:modelheavy} and \eqref{eq:gamma} hold, and $1 \le k_n \le n$, $k_n 
\to\infty$.
Then, with $\widehat \gamma_n(k)$ in \eqref{eq:def-hill-W}
\[
\begin{split}
& \lim_{n\to\infty} \frac{1}{k_n}\log P\left(\frac{\widehat\gamma_n(k_n)}{\gamma} \ge 1 +c 
\right) = -rc +r\log(1+c), \\
& \lim_{n\to\infty} \frac{1}{k_n}\log P\left(\frac{\widehat\gamma_n(k_n)}{\gamma} \le 1 -c 
\right) = rc +r\log(1-c),
\end{split}
\]
for all $c \in(0,1)$.
\end{corollary}

On the other hand, in the classical iid setup 
\cite[Theorem 2]{Cheng} proved that if the upper tail of the common 
distribution function of the underlying iid sequence is regularly varying 
with parameter $-1/\gamma$, the distribution function is continuous, and
$k_n \to \infty$, $k_n / n \to 0$ for the subsequence $k_n$, 
then
\[
\begin{split}
& \lim_{n\to\infty} \frac{1}{k_n} 
\log P\left(\frac{\widehat\gamma_n(k_n)}{\gamma} \ge 1+ c \right) =
-c +\log(1+c), \\
& \lim_{n\to\infty} \frac{1}{k_n}
\log P\left(\frac{\widehat\gamma_n(k_n)}{\gamma} \le 1- c \right) = 
c +\log(1-c)
\end{split}
\]
for all $c \in(0,1)$.
In particular, we see that large deviations do not draw distinction between 
the members in the regularly varying class, and the tail index alone determines 
the rate. On the other hand, in model \eqref{eq:modelheavy} any possible rate 
function can appear in Theorem \ref{thm:large}.
Furthermore, in the iid setting $k_n/n \to 0$ as $n \to \infty$, while in 
model \eqref{eq:modelheavy} the choice $k_n = n$ is possible, leading to 
convergence rates of $e^{-c_1 k_n}$ for some $c_1 > 0$ in the former, 
and $e^{- c_2 n}$ for some $c_2 > 0$ 
in the latter case. In Corollary \ref{cor:ldp} for $r > 1$ we 
obtain $r$ times better large deviation rates.

\subsection{Comparison of model \eqref{eq:modelheavy} with the classical iid \\
model}

Our assumptions in model \eqref{eq:modelheavy} are motivated by the theoretical result 
in Theorem \ref{thm:exponential}.
We emphasize that in the new model, we do not have a fixed
underlying distribution for the individual sample observations $W_{\delta_i, n}$,
unless $Z_1, Z_2, \ldots$ are iid exponentials. 
In the general case, the distribution of $W_{\delta_i, n}$ depends on 
$n$ and on the distribution of $Z_1$.
By Theorem \ref{thm:exponential}, as $n \to \infty$,
$W_{\delta_i, n}$ converges in distribution to the strict 
Pareto distribution with tail index $1/\gamma$. However, if $Z_1$ in 
\eqref{eq:generalized} is bounded, then $W_{\delta_i,n}$ is also bounded. Therefore,
model \eqref{eq:modelheavy} does not imply any tail asymptotics for $W_{\delta_i, n}$
for fixed $n$. 
Theorem \ref{thm:exponential} also reflects that our model is close to the iid model.

The assumptions in \eqref{eq:modelheavy} can be equivalently stated as the 
sequence of scaled log spacings $(n-k+1) (\log W_{k,n} - \log W_{k-1,n})$,
for $k=1,\ldots, n$, form
an iid sequence with finite mean. Therefore, checking whether model
\eqref{eq:modelheavy} fits the data is 
equivalent to verifying the iid assumption.
On the other hand, apart from the finiteness of 
the first moment, we do not need any assumption on the tail.

We note that we are not aware of 
any situation where assumption \eqref{eq:modelheavy} is more natural than 
the iid assumption. 
Based on Theorem \ref{thm:exponential} we expect that whenever the classical iid model 
with regularly varying tails fits the data, model \eqref{eq:modelheavy} will
fits as well. Experimental study on a few data sets also supports this statement.
We do not present new results on the Hill estimator in the iid setup. 
Instead, we propose a new model where the Hill estimator performs very well,
without any assumptions on the tail behavior.

\section{Proofs} \label{sect:proofs}

For easier reference we state a simple lemma.

\begin{lemma}
\label{lemma:permutation}
Let $(X_1, \ldots, X_n)$ be a random vector with joint density function
$f(x_1, \ldots, x_n)$.

{\rm(i)} If $\underline{v} =(v_1,\ldots,v_n)$ is a permutation of $(1,\ldots,n)$
and $(l_{1,\underline{v}},\ldots,l_{n,\underline{v}})$ is the inverse of $\underline{v}$, 
then
the joint density of $(X_{v_1}, \ldots, X_{v_n})$ is given by
\[
h(x_1,\ldots,x_n)=f(x_{l_{1,\underline{v}}}, \ldots, x_{l_{n,\underline{v}}}).
\]

{\rm(ii)} Let $(\delta_1,\ldots,\delta_n)$ denote a uniform random permutation of 
$(1,\ldots,n)$, independent of $(X_1, \ldots, X_n)$.
Then the joint density of $(X_{\delta_1}, \ldots, X_{\delta_n})$ is given by
\[
h(x_1,\ldots,x_n) =\frac{1}{n!}\sum_{\underline v}f(x_{v_1}, \ldots, x_{v_n}),
\]
where the sum $\sum_{\underline v}$ is taken over all the permutations of $(1,\ldots,n)$.
\end{lemma}

\noindent \emph{Proof.}\ 
The proof of part (i) is straightforward so that it is omitted.

For part (ii) we have for any Borel set $A$
\[
\begin{split}
& P ( (X_{\delta_1,n}, \ldots, X_{\delta_n, n}) \in A ) \\
& = 
\sum_{\underline v} 
P ( (X_{\delta_1,n}, \ldots, X_{\delta_n, n}) \in A | 
\underline{\delta} =\underline{v}) P(\underline{\delta}=\underline{v}) \\
& = \frac{1}{n!}
\sum_{\underline v}
P( (X_{v_1,n}, \ldots, X_{v_n,n}) \in A) \\
& = \frac{1}{n!} \sum_{\underline v}
\int_A f(x_{l_{1,\underline{v}}},\ldots,x_{l_{n,\underline{v}}}) d \underline x \\
& = \frac{1}{n!}\sum_{\underline v}
\int_A f(x_{v_1}, \ldots, x_{v_n}) d \underline x,
\end{split}
\]
as claimed, where the third equality follows from part (i).
\qed \medskip

For $n>1$ let $(\delta_1,\ldots,\delta_n)$ be a uniform random permutation of 
$(1,\ldots,n)$ which is independent of $(Z_1,\ldots,Z_n)$.
Let $\varphi(t) =E(\exp(itZ_1))$ denote the characteristic function of $Z_1$.
The next statement describes the finite sample properties of 
$X_{\delta_1,n},\ldots,X_{\delta_n,n}$.

\begin{lemma} \label{lemma:sample}
Assume that \eqref{eq:generalized} holds with iid random variables 
$Z_1, Z_2, \ldots$.

\begin{itemize}
\item[(i)] 
The random variables $X_{\delta_1,n},\ldots,X_{\delta_n,n}$ are identically 
distributed
and for any $1 \leq k \leq n$, the joint distribution of 
$(X_{\delta_1,n},\ldots,X_{\delta_k,n})$ is described by the characteristic function
\[
E e^{i  \sum_{j=1}^k t_j X_{\delta_j, n} }  
= \frac{(n-k)!}{n!} 
\sum_{\bl} 
\prod_{m=1}^{\max \bl} \varphi \left( \frac{\sum_{j:\ell_j \geq m} t_j}{n+1-m} \right),
\]
where summation is meant for all possible $k$-tuples $(\ell_1, \ldots, \ell_k) = \bl$ 
from $(1, \ldots, n)$ with pairwise different components.

\item[(ii)] 
If $E(|Z_1|) < \infty$, then $E(X_{\delta_1,n}) =\gamma$.

\item[(iii)] Assume that $Z_1,Z_2,\ldots$ are iid nonnegative random variables
with absolutely continuous distribution with density function $g$.
Then the distribution of $(X_{1,n},\ldots,X_{k,n})$ is absolutely continuous with density 
function
\begin{equation}
\label{eq:odensity}
p(y_1,\ldots,y_k)=\prod_{j=1}^k(n-j+1)g\big((n-j+1)(y_j-y_{j-1})\big),
\end{equation}
$0=y_0< y_1< \cdots< y_k$, and $p(y_1,\ldots,y_k)=0$ otherwise, $k=1,\ldots,n$.
The distribution of $(X_{\delta_1,n},\ldots,X_{\delta_n,n})$ is absolutely continuous with 
density function
\begin{equation}
\label{eq:density}
h(y_1,\ldots,y_n)=\prod_{j=1}^ng\big((n-j+1)(y_{j,n}-y_{j-1,n})\big),
\end{equation}
$y_1,\ldots,y_n \ge 0$, where $y_{1,n}\le\cdots\le y_{n,n}$ are the ordered values 
pertaining to $y_1,\ldots,y_n$
and $y_{0,n}=0$.
\end{itemize}
\end{lemma}

\noindent \emph{Proof.}\ 
(i)
Clearly, the random variables $X_{\delta_1,n},\ldots,X_{\delta_n,n}$ are identically 
distributed.
For $1 \leq k \leq n$, and $t_1, \ldots, t_k \in \R$ we have
\begin{equation*}
\begin{split}
E e^{i  \sum_{j=1}^k t_j X_{\delta_j, n} } & 
= \sum_{\bl} \frac{(n-k)!}{n!} 
E \exp \bigg\{ i  \sum_{j=1}^k t_j X_{\ell_j, n} \bigg\} 
\\ & 
= \frac{(n-k)!}{n!} 
\sum_{\bl} 
\prod_{m=1}^{\max \bl} \varphi \left( \frac{\sum_{j:\ell_j \geq m} t_j}{n+1-m} \right),
\end{split}
\end{equation*}
where the empty sum is defined to be 0, and summation is meant for all
$k$-tuples $(\ell_1, \ldots, \ell_k) = \bl$ from 
$(1, \ldots, n )$ with pairwise different components. 
Here we used also that
\[
\begin{split}
\sum_{j=1}^k t_j X_{\ell_j, n} &  = \sum_{j=1}^k t_j 
\sum_{m=1}^{\ell_j} \frac{Z_m}{n+1-m} 
= \sum_{m=1}^{\max \bl} \frac{Z_m}{n+1-m} \sum_{j: \ell_j \geq m } t_j.
\end{split}
\]

(ii) Again, by the law of total expectation
\[
\begin{split}
E(X_{\delta_1,n}) & =\sum_{k=1}^n E(X_{\delta_1,n}|\delta_1=k) P(\delta_1=k)
=\frac{1}{n} \sum_{k=1}^n E(X_{k,n}) \\
& =E\left(\frac{1}{n}\sum_{k=1}^n Z_k\right)
=\gamma.
\end{split}
\]

(iii) The one-to-one transformation between the vectors $(X_{1,n},\ldots,X_{k,n})$ and
$(Z_1,\ldots,Z_k)$ described in \eqref{eq:generalized} in combination with
the corresponding density transformation gives the formula in \eqref{eq:odensity}.
By Lemma \ref{lemma:permutation} (ii) we have
\[
h(y_1,\ldots,y_n) =\frac{1}{n!}\sum_{\underline v}p(y_{v_1},\ldots,y_{v_n}).
\]
This sum contains at most one nonzero component $p(y_{v_1},\ldots,y_{v_n})$ for which 
$(y_{v_1},\ldots,y_{v_n})=(y_{1,n},\ldots,y_{n,n})$. This implies \eqref{eq:density}.
\qed \medskip

\noindent \emph{Proof of Theorem \ref{thm:exponential}.}\
First let $k = 1$. Using Lemma \ref{lemma:sample} (i) with $k=1$
we have
\begin{equation} \label{eq:defpsi}
\psi_n(t) := E e^{i t X_{\delta_1, n} } = \frac{1}{n} \sum_{m=1}^n 
\prod_{j=1}^m \varphi\left( \frac{t}{n+1 -j} \right).
\end{equation}
We prove the convergence of the characteristic function.
Fix $\varepsilon > 0$ arbitrary small, and $t \in \R$. 
Since $\varphi(s) = 1 + i  s \gamma + o(s)$ as $s \to 0$,
for $j \leq (1 - \varepsilon ) n$ 
\begin{equation} \label{eq:phi-asy}
\log \varphi \left( \frac{t}{n+1 - j} \right) = 
i  t \gamma (n+1 - j)^{-1} + r(j,n),
\end{equation}
where for $n$ large enough
$|r(j,n)| \leq \varepsilon^2 |t| (n+1-j)^{-1} \leq \varepsilon |t| n^{-1}$.
Therefore, for $m \leq (1- \varepsilon) n$
\begin{equation*}
\begin{split}
\prod_{j=1}^m  \varphi \left( \frac{t}{n+1 - j} \right)
& = \exp \left\{ 
\sum_{j=1}^m \log \varphi \left( \frac{t}{n+1 - j} \right)
\right\} \\
& = \exp \left\{ i  t \gamma \sum_{j=1}^m (n+1 - j)^{-1} + R(n,m) \right\} \\
& = \exp \left\{ - i  t \gamma \log \left( 1 - \frac{m}{n} \right)  
+ R_2(n,m) \right\},
\end{split}
\end{equation*} 
where by \eqref{eq:phi-asy} $R(n,m) \leq \varepsilon |t|$. Furthermore,
for $m \leq (1-\varepsilon) n$
\begin{equation} \label{eq:log-asy}
 \int_{n-m}^n x^{-1} d x - \sum_{j={n-m+1}}^n j^{-1} \leq \frac{1}{n \varepsilon},
\end{equation}
implying that $|R_2(n,m)| \leq 2 \varepsilon |t|$ for $n$ large enough.
Thus
\begin{equation} \label{eq:char-fgv-bound1}
\begin{split}
& \left| \sum_{m=1}^{(1-\varepsilon)n} 
\left[ \prod_{j=1}^m \varphi \left( \frac{t}{n+1 - j} \right) -
\exp \left\{ - i  t \gamma \log \left( 1 - \frac{m}{n} \right) \right\} 
\right] \right| \\
& \leq \sum_{m=1}^{(1-\varepsilon)n} \left| e^{R_2(n,m)} - 1 \right| 
\leq 2 \varepsilon |t| n.
\end{split}
\end{equation}

Next, by Riemann approximation we get as $n \to \infty$
\begin{equation} \label{eq:char-fgv-conv}
\begin{split}
\frac{1}{n} \sum_{m=1}^{(1-\varepsilon) n} 
\exp \left\{ - i  t \gamma \log \left( 1 - \frac{m}{n} \right) \right\} 
& \to \int_\varepsilon^1 e^{- i  t \gamma \log x } d x \\
& = \int_0^{- \gamma \log \varepsilon}  
e^{i  t y} e^{-y/\gamma} \gamma^{-1} d y.
\end{split}
\end{equation}
As $\varepsilon \to 0$ the last integral converges to 
$(1 - i  t \gamma)^{-1}$, the characteristic function of the 
exponential distribution with mean $\gamma$.

Recalling \eqref{eq:defpsi}, we may write
\[
\begin{split}
& \psi_n(t) - \frac{1}{1 - i  t \gamma} \\
& = 
\frac{1}{n}\sum_{m=1}^{(1-\varepsilon)n} 
\prod_{j=1}^m \varphi \left( \frac{t}{n+1 - j} \right) 
- \int_0^{- \gamma \log \varepsilon}  
e^{i  t y} e^{-y/\gamma} \gamma^{-1} d y \\
& \quad + \frac{1}{n}\sum_{m=(1-\varepsilon)n}^n 
\prod_{j=1}^m \varphi \left( \frac{t}{n+1 - j} \right) 
 - \int_{- \gamma \log \varepsilon}^{\infty} 
e^{i  t y} e^{-y/\gamma} \gamma^{-1} d y.
\end{split}
\]
We see that the difference of the first two terms is 
bounded by $3 \varepsilon |t|$ for $n$ large enough
by \eqref{eq:char-fgv-bound1} and  \eqref{eq:char-fgv-conv}, 
while the third and fourth term tend to 0 as $\varepsilon \downarrow 0$.
Thus the convergence holds for $k = 1$.
\medskip

For $k > 1$ the proof is essentially the same.  Fix $t_1, \ldots, t_k \in \R$,
and $\varepsilon > 0$ small.
To ease notation put $T(\bl, m) = 
\sum_{j: \ell_j \geq m } t_j / (n+1-m)$. 
Then, similarly as in \eqref{eq:phi-asy}
for $m \leq (1- \varepsilon) n$,
\[
\log \varphi( T(\bl, m)) = i  \gamma T(\bl,m) + r(\bl, m),
\]
with $|r(\bl, m)| \leq \varepsilon / n$. Thus
\[
\sum_{m=1}^{(1-\varepsilon ) n} 
\log \varphi( T(\bl, m)) = i  \gamma \sum_{m=1}^{(1-\varepsilon ) n} T(\bl,m) 
+ R(\bl),
\]
where $|R(\bl)| \le \varepsilon$ for $n$ large enough.
It is easy to see that
\[
\sum_{m=1}^{\max \bl} T(\bl,m) = \sum_{j=1}^k t_j \sum_{m=1}^{\ell_j} \frac{1}{n+1-m}.
\]
Assuming $\max_j \ell_j \leq (1- \varepsilon) n$, using \eqref{eq:log-asy}
we obtain as in \eqref{eq:char-fgv-bound1}
\[
\begin{split}
& \left| \sum_{\bl: \max \bl \leq (1- \varepsilon)n} 
\left[
\prod_{m=1}^{\max \bl} \varphi( T(\bl, m) ) - 
\exp \bigg\{ - i  \gamma 
\sum_{j=1}^k t_j \log \left( 1 - \frac{\ell_j}{n} \right) \bigg\} \right] \right| \\
& \leq  \varepsilon n^k  (k+1) \sum_{j=1}^k |t_j| .
\end{split}
\]
Furthermore,
\[
\begin{split}
& 
\frac{1}{n (n-1) \ldots (n-k+1)}
\sum_{\bl: \max \bl \leq (1- \varepsilon)n} 
\exp \bigg\{ - i  \gamma 
\sum_{j=1}^k t_j \log \left( 1 - \frac{\ell_j}{n} \right) \bigg\} \\
& \sim 
\prod_{j=1}^k 
\frac{1}{n}
\sum_{\ell=1}^{(1- \varepsilon)n} 
\exp \bigg\{ - i  \gamma t_j \log \left( 1 - \frac{\ell}{n} \right) \bigg\} \\
& \to \prod_{j=1}^k 
\int_0^{- \gamma \log \varepsilon}  
e^{i  t_j y} e^{-y/\gamma} \gamma^{-1} d y,
\end{split}
\]
where the last convergence holds as $n \to \infty$, and
$\sim$ stands for asymptotic equality, that 
is $a_n \sim b_n$ means that $\lim_{n \to \infty} \tfrac{a_n}{b_n} = 1$.
As $\varepsilon \downarrow 0$ the right-hand side converges to the 
characteristic function of independent exponentials with mean $\gamma$, the 
proof can be finished as in the $k =1$ case.
\qed \medskip

\noindent \emph{Proof of Theorem \ref{thm:moments}.}\ 
Set $m_{k,n} =E\big(X_{\delta_1,n}^k\big)$ and $\mu_k=E\big(Z_1^k\big)$.
By Lemma \ref{lemma:sample} (ii), $m_{1,n} =\gamma$.
By induction, we will show that $m_{k,n}\to k!\gamma^k$.
By \eqref{eq:defpsi} we have
\[
\begin{split}
\psi_n(t) =
\frac{1}{n}\varphi\left(\frac{t}{n}\right) 
+\frac{n-1}{n}\varphi\left(\frac{t}{n}\right)\psi_{n-1}(t).
\end{split}
\]
The Leibniz rule for the derivatives of a product function gives that the $k$th derivative of $\psi_n(t)$ is
\[
\psi_n^{(k)}(t) =\frac{1}{n^{k+1}}\varphi^{(k)}\left(\frac{t}{n}\right)
+\frac{n-1}{n}\sum_{j=0}^k\binom{k}{j}\frac{1}{n^{k-j}}\varphi^{(k-j)}\left(\frac{t}{n}\right)\psi_{n-1}^{(j)}(t).
\]
Then
\[
m_{k,n}=\frac{\psi_n^{(k)}(0)}{i^k} =\frac{\mu_k}{n^k}
+\frac{n-1}{n}m_{k,n-1} +\frac{n-1}{n}\sum_{j=1}^{k-1}\binom{k}{j}\frac{1}{n^{k-j}}\mu_{k-j}m_{j,n-1}.
\]
For $k>1$ this recursion gives
\begin{equation}
\label{eq:recurs}
\begin{split}
m_{k,n} =&
\frac{\gamma k}{n}\sum_{r=1}^{n-1} \frac{n-r}{n-r+1} m_{k-1,n-r}
+\frac{\mu_k}{n} \sum_{r=0}^{n-1}\frac{1}{(n-r)^{k-1}} \\
&+\sum_{j=1}^{k-2}\sum_{r=0}^{n-2}\frac{n-r-1}{n}\binom{k}{j}\frac{1}{(n-r)^{k-j}}\mu_{k-j}m_{j,n-r-1}.
\end{split}
\end{equation}
Assuming $m_{k-1,n}\to (k-1)!\gamma^{k-1}$, the first term in the right side of
\eqref{eq:recurs} converges to $k!\gamma^k$,
and the second term tends to zero.
Further, for a fixed $j$,
$1\le j\le k-2$, by the induction hypothesis we may assume that $m_{j,n}\to j!\gamma^j$.
Hence $|m_{j,n}|\le M_j$ for some bound $M_j$. It follows that
\[
\begin{split}
\Bigg|\sum_{r=0}^{n-2}\frac{n-r-1}{n}\binom{k}{j}\frac{1}{(n-r)^{k-j}}&\mu_{k-j}m_{j,n-r-1}\Bigg|\\
\le M_j &\mu_{k-j}\binom{k}{j}\sum_{r=0}^{n-2}\frac{n-r-1}{n}\frac{1}{(n-r)^{k-j}}.
\end{split}
\]
The last term and hence the double sum in \eqref{eq:recurs} tends to zero.
This imply the statement $m_{k,n}\to k!\gamma^k$.

To see the asymptotic correlation structure, it suffices to prove 
the convergence 
$E(X_{\delta_1,n}X_{\delta_2,n}) \to\gamma^2$.
Write
\[
\begin{split}
& E(X_{\delta_1,n}X_{\delta_2,n}) =
E\left( \left(\frac{Z_1}{n}+\widetilde 
X_{\delta_1-1,n-1}\right) \left(\frac{Z_1}{n}+\widetilde X_{\delta_2-1,n-1}\right) 
\right)\\
& =E\left(\frac{Z_1^2}{n^2}\right) +2E\left(\frac{Z_1}{n}\right)E\left(\widetilde X_{\delta_1-1,n-1}\right)
+ E\left(\widetilde X_{\delta_1-1,n-1}\widetilde X_{\delta_2-1,n-1}\right),
\end{split}
\]
where
\begin{equation*}
\widetilde X_{\delta_\ell-1,n-1} = X_{\delta_\ell,n} - \frac{Z_1}{n}
=\sum_{j=2}^{\delta_\ell}\frac{Z_j}{n+1-j}, \quad \ell = 1,2.
\end{equation*}
Note that given $\delta_1>1$, 
$\widetilde X_{\delta_1-1,n-1} \stackrel{\mathcal{D}}{=} 
X_{\delta_1,n-1}$, and given $\delta_1 > 1, \delta_2 > 1$,
$(\widetilde X_{\delta_1 -1, n-1}, \widetilde X_{\delta_{2} - 1, n-1}) 
\stackrel{\mathcal{D}}{=} (X_{\delta_1, n-1}, X_{\delta_2, n-1})$.
Therefore,  $E(\widetilde X_{\delta_k-1,n-1}) =$ $E\left(X_{\delta_k,n} 
-\frac{Z_1}{n}\right)=\gamma-\frac{\gamma}{n}$
and
\begin{align}
E\big(\widetilde X_{\delta_1-1,n-1}\widetilde X_{\delta_2-1,n-1}\big)
= &E\big(\widetilde X_{\delta_1-1,n-1}\widetilde X_{\delta_2-1,n-1}|\delta_1\!\ne\!1,\delta_2\!\ne\!1\big)\frac{n-2}{n}\nonumber\\
\label{eq:cov}
& +E\big(\widetilde X_{\delta_1-1,n-1}\widetilde X_{\delta_2-1,n-1}|\delta_1\!=\!1\ {\rm or}\ \delta_2 =1\big)\frac{2}{n}.
\end{align}
The quantity in \eqref{eq:cov} equals zero.
Then, writing $C_n =E(X_{\delta_1,n}X_{\delta_2,n})$,
\[
C_n =\frac{\mu_2}{n^2} +2\frac{\gamma}{n}\left(\gamma-\frac{\gamma}{n}\right)+C_{n-1}\frac{n-2}{n}.
\]
This recursion yields
\[
C_n =\frac{2C_2}{n(n-1)} +\sum_{j=0}^{n-3}\frac{\mu_2+2\gamma^2(n-j-1)}{(n-j)}\frac{(n-j-1)}{n(n-1)},
\]
where $C_2 =E(X_{1,2}X_{2,2}) =\frac{\mu_2}{4}+\frac{\gamma^2}{2}$.
A routine calculation gives that $C_n\to\gamma^2$.
\qed \medskip

\noindent \emph{Proof of Theorem \ref{thm:Qn}.}\ 
By definition \eqref{eq:quant}
\begin{equation} \label{eq:vars}
\var ( Q_n^{(1)}(s) ) = \sum_{j=1}^{\lceil ns \rceil} 
\frac{\sigma^2}{(n+1-j)^2} \sim \frac{\sigma^2}{n} \frac{s}{1-s}.
\end{equation}
We show that  Lindeberg's condition holds
for $\sum_{j=1}^{\lceil ns \rceil} \tfrac{Z_j - \gamma}{n+1-j}$.
Fix $\varepsilon > 0$. Putting $s_n^2 = \var(Q_n^{(1)}(s))$, and $\ind$
for the indicator of an event, we have
\[
\begin{split}
& \frac{1}{s_n^2} \sum_{j=1}^{\lceil ns \rceil} 
E \left( \frac{(Z_j - \gamma)^2}{(n+1-j)^2}
\ind( |Z_j - \gamma| > \varepsilon s_n (n+1 - j) \right) \\
& \leq \frac{\lceil ns \rceil }{s_n^2 (n + 1 - \lceil ns \rceil)^2 }
E \left( (Z_1 - \gamma)^2 
\ind ( |Z_1 - \gamma| > \varepsilon s_n (n + 1 - \lceil ns \rceil)) \right) \to 0,
\end{split}
\]
where we used that by \eqref{eq:vars}, $s_n \sim \sigma \sqrt{s} / \sqrt{n (1-s)}$.
Therefore, by Lindeberg's central limit theorem
\[
\frac{1}{s_n} \sum_{j=1}^{\lceil ns \rceil} \frac{Z_j-\gamma}{n+1-j} 
\stackrel{\mathcal{D}}{\longrightarrow} N(0, 1).
\]
Noting that 
\[
\lim_{n \to \infty} 
\sqrt{n} \left( \sum_{j=1}^{\lceil ns \rceil} \frac{1}{n+1 -j} + \log ( 1- s) \right) = 0,
\]
we proved that 
\[
\sqrt{n} \frac{- \sqrt{1-s} \log ( 1 - s)}{ \sqrt{s} \sigma}
( \widetilde \gamma_n(s) - \gamma ) \stackrel{\mathcal{D}}{\longrightarrow} N(0, 1),
\]
as claimed.
\qed \medskip

Proposition \ref{prop:likelihood} is a straightforward consequence of the following.

\begin{lemma} \label{lemma:cond}
Assume that \eqref{eq:generalized} holds, and the random variables $Z_1,Z_2,\ldots$ 
have density $g$. Then the conditional distribution of 
$(X_{n-k+1,n},\ldots,X_{n,n})$ given $X_{n-k,n} =x_{n-k}$ is absolutely continuous with 
density function
\begin{equation}
\widetilde h(x_{n-k+1},x_{n-k+2},\ldots,x_n|x_{n-k})
=k!\prod_{j=n-k+1}^ng\big((n-j+1)(x_j-x_{j-1})\big),
\end{equation}
for $0\le x_{n-k}\le x_{n-k+1}\le\cdots\le x_n$, $k=1,\ldots,n$, $x_0=0$.
\end{lemma}

\noindent \emph{Proof.}\ 
Based on representation \eqref{eq:generalized}, it is easy to check the Markov 
property of the statistics $X_{1,n},\ldots,X_{n,n}$.
For all $y_1,\ldots,y_n\in\mathbb R$, $k\in \{1, \ldots, n\}$,
\begin{equation}
\label{eq:Markov}
\begin{split}
P(&X_{k,n}<y_k,\ldots,X_{n,n}<y_n|X_{1,n}=y_1,\ldots,X_{k-1,n}=y_{k-1})\\
& =P(X_{k,n}<y_k,\ldots,X_{n,n}<y_n|X_{k-1,n}=y_{k-1}).
\end{split}
\end{equation}

Let $h(x_{i_1},\ldots,x_{i_l}|x_{j_1},\ldots,x_{j_s})$,
$1 \leq j_1 < \cdots < j_s < i_1 < \cdots < i_l \leq n$
denote the conditional density of $(X_{i_1,n},\ldots,X_{i_l,n})$ given 
$X_{j_1,n}=x_{j_1},\ldots$, $X_{j_s,n}=x_{j_s}$. Then using the Markov property 
\eqref{eq:Markov} and the density formula in \eqref{eq:odensity},
\[
\begin{split}
h &(x_{n-k+1},x_{n-k+2},\ldots,x_{n}|x_{n-k})\\
& =h(x_{n-k+1},x_{n-k+2},\ldots,x_{n}|x_{1},\ldots,x_{n-k})
=\frac{p(x_{1},\ldots,x_{n})}{p(x_{1},\ldots,x_{n-k})}\\
& =k!\prod_{j=n-k+1}^n g\big((n-j+1)(x_{j}-x_{j-1})\big).
\end{split}
\]
\qed 
\bigskip
\bigskip

\noindent \textbf{Acknowledgements}

\bigskip

We are grateful to the anonymous referees and to the AE for helpful 
comments and suggestions, which greatly improved the presentation of our results.

This research was supported by the Ministry of Culture and Innovation of 
Hungary from the National Research, Development
and Innovation Fund, project no.~TKP2021-NVA-09. PK was supported by the
J\'{a}nos Bolyai Research Scholarship of the Hungarian Academy of Sciences.

The R codes are publicly available at 
https://github.com/pkevei/Renyi-simulation/blob/main/Renyi-code.r.

\end{document}